 \documentclass[11pt]{article}
 \usepackage[top=1in,bottom=1in,left=1.5in,right=1.5in]{geometry}
  \usepackage{amsmath,amssymb}
  \usepackage{latexsym}% defines $\Box$ for LaTeX2e
  \usepackage{graphicx,color}

  \newcommand{\C}{\mathbb{C}}
  \newcommand{\F}{\mathbb{F}}
  \newcommand{\N}{\mathbb{N}}
  \renewcommand{\P}{\mathbb{P}}
  \newcommand{\R}{\mathbb{R}}
  \newcommand{\Z}{\mathbb{Z}}

  \newcommand{\e}{\mathbf{e}}
  \newcommand{\f}{\mathbf{f}}
  \newcommand{\g}{\mathbf{g}}
  
  \newcommand{\m}{\mathbf{m}}

  \newcommand{\U}{\mathbf{U}}
  \newcommand{\uu}{\mathbf{u}}
  \newcommand{\vv}{\mathbf{v}}
  
  \newcommand{\w}{\mathbf{w}}
  \newcommand{\W}{\mathbf{W}}
  \newcommand{\x}{\mathbf{x}}
  
  \newcommand{\y}{\mathbf{y}}
  \newcommand{\z}{\mathbf{z}}
  \newcommand{\0}{\mathbf{0}}

  \newcommand{\cA}{\mathcal{A}}
  \newcommand{\cB}{\mathcal{B}}
  \newcommand{\cC}{\mathcal{C}}

  \newcommand{\cF}{\mathcal{F}}
  
  \newcommand{\cG}{\mathcal{G}}

  \newcommand{\cO}{\mathcal{O}}

  \newcommand{\cR}{\mathcal{R}}
  \newcommand{\cS}{\mathcal{S}}
  \newcommand{\cT}{\mathcal{T}}

  \newcommand{\cX}{\mathcal{X}}

  \newcommand{\hs}{\hspace*{\parindent}}
  \newcommand{\proof}{\hs \textbf{Proof.\ }}
  
  \newcommand{\tr}{\mathop{\mathrm{tr}}\nolimits}

  \newcommand{\trans}{^\top}
  \newcommand{\qed}{\hspace*{\fill} $\Box$\\}

  \newcommand{\dist}{\mathrm{dist}}

  \newcommand{\Krank}{\mathrm{krank\;}}

  \newcommand{\range}{\mathrm{range\;}}
  \newcommand{\Sing}{\mathrm{Sing\;}}

  \newcommand{\rB}{\mathrm{B}}

  \newcommand{\rH}{\mathrm{H}}

  \newcommand{\rS}{\mathrm{S}}

  \newcommand{\rank}{\mathrm{rank\;}}
  \newcommand{\srank}{\mathrm{srank\;}}

  \newtheorem{theo}{\bfseries \hs Theorem}[section]
  \newtheorem{defn}[theo]{\bfseries \hs Definition}
  \newtheorem{prop}[theo]{\bfseries \hs Proposition}
  \newtheorem{lemma}[theo]{\bfseries \hs Lemma}
  
  \newtheorem{corol}[theo]{\bfseries \hs Corollary}
  
  \newtheorem{example}[theo]{\bfseries \hs Example}

  \numberwithin{equation}{section} % Automatically number equations within sections

 \renewcommand{\span}{\mathrm{span}}

\begin{document}

 \title{Best approximation on semi-algebraic sets and\\$k$-border rank approximation of symmetric tensors}

 \author
 {Shmuel Friedland\footnote{Department of Mathematics, Statistics and Computer Science,
 University of Illinois at Chicago,
 Chicago, Illinois 60607-7045, USA,
 \emph{email}: friedlan@uic.edu. This work was supported by NSF grant DMS-1216393} \ and Ma\l gorzata Stawiska
 \footnote{Mathematical Reviews, 416 Fourth Street, Ann Arbor, MI 48103-4816, USA, \emph{email}: stawiska@umich.edu}}

 \date{November 6, 2013}
 \maketitle

 \begin{abstract}
 In the first part of this paper we study a best approximation of a vector in Euclidean space $\R^n$ with respect to a closed semi-algebraic set $C$
 and a given semi-algebraic norm.  Assuming that the given norm and its dual norm are differentiable we show that 
 a best approximation is unique outside a hypersurface.  We then study the case where $C$ is an irreducible variety and the approximation is with respect  
to the  Euclidean norm.
We show that for a general point in $\x\in \R^n$  the number of critical points of the distance function of $\x$ to $C$ is bounded above by a degree 
of a related dominant map.  If $C$ induces a smooth projective variety $ V_\P\subset P(\C^n)$ then this degree is the top Chern number of a corresponding
vector bundle on $V_\P$.  We then study the problem when a best $k(\ge 2)$-border rank approximation of a symmetric tensor is symmetric.  We show that under certain
dimensional conditions there exists an open semi-algebraic set of symmetric tensors for which a best $k$-border rank is unique and symmetric.

 \end{abstract}



 \section{Introduction}\label{sec:intro}

 Let $\nu:\R^n\to [0,\infty)$ be a norm on $\R^n$.
 In many applications one needs to approximate a given vector $\x\in\R^n$ by a point $\y$ in a given closed subset $C\subset \R^n$.
 Assume that we measure the approximation of $\y$ to $\x$ by $\nu(\x -\y)$.  Then the distance of $\x$ to $C$ with respect to the norm
 $\nu$ is defined as $\dist_{\nu}(\x,C):=\min \{\nu(\x-\y), \;\y\in C\}$.   A point $\y^\star\in C$ is called a best $\nu$-approximation of $\x$ 
 if $\nu(\x-\y^{\star})=\dist_{\nu}(\x,C)$.  Let $\|\cdot\|$denote  the Euclidean norm on $\R^n$ and let   $\dist(\x,C)$ denote
 the distance $\dist_{\|\cdot\|}(\x,C)$.  We call a best $\|\cdot\|$-approximation  a best $C$-approximation, or briefly  a best approximation. \\
 
 One of important applications is an approximation of 
 a given matrix $A\in\R^{p\times q}$ by a matrix of rank at most $k$ in the Frobenius norm., i.e., the Euclidean norm on $\R^{p\times q}$.  
 This classical problem has a well-known solution, namely, the singular value decomposition (SVD) of $A$ gives a best rank $k$ 
 approximation \cite{GV}.  Recall that finding SVD decomposition with in $\varepsilon$ precision has a polynomial time algorithm.
On the other hand, a similar problem for tensors, i.e., approximating a given $d$-mode tensor for $d\ge 3$, is much more difficult
\cite{HL13}.\\

In general, the numerical methods for finding best approximation to $\x$ in $C$ are based on finding a local minimum of the function $\nu(\x-\y), \y\in C$.
For example, to find a best rank one approximation of tensors one uses an alternating least squares method \cite{FMPS}.
Most of these methods at most will converge to a local minimum point of $f_{\x,\nu}(\y):=\nu(\x-\y), \y\in C$.  
We use the abbreviation $f_{\x}$ for the function $f_{\x,\|\cdot\|}$.

In a recent paper,  the first named author and G. Ottaviani \cite{FO12}  considered the above approximation problem  in the following special setting. 
Let $\m=(m_1,\ldots,m_d)\in\N^d$.  For a field $\F$ let $\F^{\m}\equiv\F^{m_1\times\ldots\times m_d}:=\otimes_{i=1}^d \F^{m_i}$.  Let
$\cR(1,\m),\cR_\C(1,\m)$ be the closed sets of all decomposable tensors $\otimes_{i=1}^d \x_i$ in  $\otimes_{i=1}^d\R^{m_i},
\otimes_{i=1}^d\C^{m_i}$ respectively.   For $C=\cR(1,\m)$ a best approximation to $\cT\in\R^{\m}$ is called a best rank one approximation.
Critical points of $f_{\cT}(\cX):=\|\cT-\cX\|, \cX\in\cR(1,\m)$ were characterized by L.-H. Lim \cite{Lim05}, which gave rise to the definition of
singular vectors and singular values for tensors.  The notion of singular vectors can be naturally extended to complex-valued tensors $\cT\in\C^\m$.
It was shown in \cite{FO12} that  there exists a variety $U\subset\C^\m$ such that each $\cT\in\C^\m\setminus U$
has exactly $\delta(\cR_\C(1,\m))$ critical points.
Furthermore, a closed formula for $\delta(\cR_\C(1,\m))$ is given in \cite{FO12}.  Moreover,
a best rank approximation is unique for almost all $\cT\in\R^{\m}$.  \\

Let $m^{\times d}:=(\underbrace{m,\ldots,m}_d)$ and let 
$\rS(d,m,\R)\subset \R^{m^{\times d}}$ denote  the $d$-mode symmetric tensors.  
Recall Banach's theorem \cite{Ban38}, which can be stated as follows:  
Every symmetric tensor  $\cS\in\rS(d,m,\R)$ has a symmetric best rank one approximation.
The Banach theorem was re-proved in \cite{CHLZ, Fri13}.
It is shown in \cite{FO12} that for almost all $\cS\in\rS(m,d,\R)$
a best rank one approximation is unique and symmetric.\\

Let $\cR(k,\m)\subset\R^\m$ be the closure of all tensors of rank at most $k$, i.e., all tensors of border rank at most $k$.
For $C=\cR(k,\m)$ we call a best approximation by a best $k$-border rank approximation.
It is natural to ask if Banach's theorem holds for $C=\cR(k,m^{\times d}), k\ge 2$.  I.e., does there exist a symmetric best $k$-border
rank approximation to every symmetric tensor $S\in\rS(d,m,\R)$?  We give a partial positive answer to this problem for certain values
of $k$.  Namely, there exists function $N(m,d)\in\N$, given by \eqref{defNnd}, such that for each $k\le N(d,m)$ 
there exists an open semi-algebraic set in $\rS(d,m,\R)$ 
for which a best $k$-border rank approximation is unique and symmetric.\\

We now summarize  the contents of the paper.  In \S\ref{sec:uniqbnuap} we discuss the best $\nu$-approximation to a closed subset $C\subset \R^n$.
We show that $\dist_\nu(\cdot,C)$ is Lipschitz on $\R^n$. 
Assuming that $\nu$  and its dual norm $\nu^*$ are differentiable we show that for a.a. $\x\in \R^n$ the $\nu$ approximation of $\y^\star\in C$ is unique.
In \S\ref{sec:semalg} we assume that $C$ is a semi-algebraic set and the norm $\nu$ is a semi-algebraic function.  We show that the function $\dist_\nu(\cdot,C)$ is 
a semi-algebraic function.  Furthermore, if $\nu$ and $\nu^*$ are differentiable then the set of $\x\in\R^n$ which do not have a unique best $\nu$ approximation 
is a semi-algebraic set of dimension less than $n$.  In \S\ref{sec:irvar} we discuss the case where $C$ is an irreducible algebraic variety in $\R^n$ and $\nu(\cdot)=
\|\cdot\|$.  For each $\x\in\R^n$ we characterize the critical points of the function $f_\x$ on  $C$.  The notion of critical points of $f_\x$ extends naturally
to $\z\in\C^n$ with respect to $C_\C\subset\C^n$, the corresponding complex irreducible variety in $\C^n$.  Assume for simplicity of the exposition that $C_\C$ is
smooth.  Let $\Sigma(C_\C)\subset \C^n\times C_\C$ be the variety of the points $(\z,\y)$ where $\y$ is a critical point of $f_\z$.   Then $\Sigma(C_\C)$
is an irreducible variety of dimension $n$.  Furthermore, the projection of $\Sigma(C_\C)$ on its first $n$ coordinates is a dominant map.  Let $\delta(C_\C)$
be the degree of this map.  Hence there exists a subvariety $U\subset \C^n$ such that for each $\z\in\C^n\setminus U$, $f_\z$ has exactly $\delta(C_\C)$ critical
points.  In particular, for each $\x\in\R^n\setminus U$ the real function $f_\x$ has at most $\delta(C_\C)$ critical points in $C$.  The case where $C_\C$ has singular
points is also discussed.  In \S\ref{sec:smtprojvar} we assume that $C_\C$ is a homogeneous variety such that the projective variety $(C_\C)_\P$ is smooth.
Suppose furthermore that $(C_\C)_\P$ intersects the projective standard quadric 
$Q_n:=\{\sum_{i=1}^n x_i^2=0, \x=(x_1,\ldots,x_n)\trans\in \C^n\}$ transversally.
Then we define a vector bundle $\cC$ over $(C_\C)_\P$ of rank $\dim (C_\C)_\P$.  We show that the top Chern number of $\cC$ is $\delta(V)$.
In \S\ref{sec:distsubs} we consider a nontrivial  finite group of orthogonal matrices acting as linear transformations on $\R^n$ and  
keeping the irreducible variety $C$  fixed.  Let $\U\subset \C^n$ be the subspace of a the fixed points of $\cG$.  Assume that $\dim \U\in\{1,\ldots,n-1\}$.
Let $\U_\R:=\U\cap\R^n$.  (The complex dimension of $\U$ is equal to the real dimension of $\U_\R$.)  We study the question when 
there exists a best $C$-approximation of $\x\in \U_\R\setminus C$ in $C\cap\U_\R$.  (This problem is a generalization of the problem when
a best $k$-border rank approximation of a symmetric tensor $\cS$ is symmetric.)
 We show that under the assumption that $C\cap\U_\R$ contains smooth points
of $C$ there exists an open semi-algebraic set $O\subset \R^n$ of dimension $n$ containing all smooth point of $C$ lying in $C\cap\U_\R$ such that for each
$\x\in O\cap\U_\R$ the best $C$ approximation of $\x$ is unique and lies in $C\cap\U_\R$.  In \S\ref{sec:genrankt} we first discuss the notion of symmetric
tensor rank of symmetric tensor over a field $\F$.  Then we define the notion of a generic tensor in $\F^{m^{\times d}}$ given as a sum
of $r$ rank one tensors, and of a generic symmetric tensor $\rS(m,d,\F)$ given as a linear combination of $r$ symmetric rank one tensors. 
Using Kruskal's theorem we show that for certain values of $r$ the generic tensors have rank $r$.  In \S\ref{sec:aprrten}
we show that for an integer $k, 2\le k\le N(m,d),$ there exists an open semi-algebraic set $ O\subset\rS(m,d,\R)$ containing all generic symmetric tensors of rank $k$
such that a best k-border rank approximation for each $\cS\in O$ is unique and symmetric.

After we completed this paper we became aware of \cite{CF00, DHOST} which are related to some of the results of our paper.

\section{Uniqueness of a best $\nu$-approximation}\label{sec:uniqbnuap}

Let $\nu$ be a norm on $\R^n$.  Let
\begin{equation}\label{defBSnu}
\rB_{\nu}:=\{\x\in \R^n, \nu(\x)\le 1\}, \quad \rS_{\nu}:=\{\x\in\R^n,\nu(\x)=1\},
\end{equation} 
denote respectively the unit ball and the unit sphere with respect to $\nu$.  It is well known that all norms on $\R^n$ are equivalent, i.e.,

\begin{equation}\label{equivnorm}
\kappa_1(\nu)\|\x\|\le \nu(\x)\le \kappa_2(\nu)\|\x\| \textrm{ for all } \x\in\R^n, \textrm{ where } 0<\kappa_1(\nu)\le \kappa_2(\nu).
\end{equation}

Recall that the dual norm $\nu^*$ is defined as $\nu^*(\x)=\max_{\y\in\rS_{\nu}} \y\trans\x$.
Since $\nu$ is a convex function on $\R^n$ it follows that  the hyperplane $\y\trans \z=1, \z\in\R^n$ is a supporting hyperplane of $B_\nu$ at $\x\in\rS_\nu$
if and only if $\y\trans \x=1$ and $\y\in \rS_{\nu^*}$.  Furthermore, $\nu$ is differentiable at $\x$ if and only if the supporting hyperplane at $\x$ is unique
\cite{Roc}.  We say that $\nu$ is differentiable if it is differentiable on $\R^n\setminus\{\0\}$.  Since $\nu$ is homogeneous, $\nu$ is 
differentiable if and only if at each $\x\in \rS_\nu$ the supporting hyperplane is unique. Assume that $\nu$ is differentiable.  For $\x\ne \0$ denote by
$\partial\nu(\x)$ the differential of $\nu$ at $\x$.  We view $\partial\nu(\x)$ as a row vector in $\R^n$.  So the directional derivative of $\nu$ at $\x$ in the 
direction $\uu\in\R^n$ is given as $ \partial\nu(\x)\uu$.
 Note that the $\ell_p$-norm on $\R^n$, 
$\|(x_1,\ldots,x_n)\trans\|_p=(\sum_{i=1}^n |x_i|^p)^{\frac{1}{p}}$, is differentiable for $p\in (1,\infty)$.\\

 The following  are generalizations of the results in \cite[\S6]{FO12}.

 \begin{theo}\label{bestnraprxlemma}  Let $\nu$ be a differentiable norm on $\R^n$ and $C\subsetneq\R^n$ be a nonempty closed set.
Let $\U\subset \R^n$ be a subspace with $\dim \U\in \{1,...,n\}$ and such that $\U$ is not contained in $C$.
Let $d(\x), \x\in \U$ be the restriction of $\dist_{\nu}(\cdot,C)$ to $\U$.   Then
\begin{enumerate}
\item%{(i)} 
\begin{equation}\label{Lipscon}
|\dist_{\nu}(\x,C)-\dist_{\nu}(\z,C)|\le \nu(\x-\z) \textrm{ for all } \x,\z\in \R^n.
\end{equation}
Hence $\dist_{\nu}(\x,C)$ is Lipschitz.
\item%{(ii)} 
The function $d(\cdot)$ is differentiable a.e. in $\U$.
\item%{(iii)} 
Let $\x\in\U\setminus C$ and assume that $d(\cdot)$ is differentiable at $\x$.
Let $\partial d(\x)$  denote the differential at $\x$ which is viewed as linear functional
 on $\U$. Let $\y^\star \in C $ be a best $\nu$-approximation to $\x$.  Then
 \begin{equation}\label{forderdistV}
 \partial d(\x)(\uu)= \partial\nu(\x-\y^\star)\uu \textrm{ for each } \uu\in\U.
 \end{equation}
 If $\z^{\star}$ is another best $\nu$-approximation to $\x$ then $\partial\nu(\x-\z^\star)-\partial\nu(\x-\y^\star)$ is orthogonal to $\U$.

\end{enumerate}

Suppose furthermore that $\nu^*$ is differentiable.  Then  at each point where $\dist_{\nu}(\x,C)$ is differentiable a best $\nu$-approximation 
is unique.  

%In particular, for a.a. $\x\in\R^n$ a best $\nu$-approximation is unique 

\end{theo}
\proof 
Let $\y^\star$ be a best $\nu$-approximation to $\x\in\R^n$.  Then 
\begin{eqnarray*}
&&\dist_\nu(\z,C)\le \nu(\z-\y^\star)= \nu(\z-\x+\x-\y^\star)\le  \nu(\z-\x)+\nu(\x-\y^\star)=\\
&&\nu(\z-\x)+\dist_\nu(\x,C) \Rightarrow \dist_\nu(\z,C)-\dist_\nu(\x,C)\le \nu(\z-\x)\Rightarrow\\ 
&&\dist_\nu(\x,C)-\dist_\nu(\z,C)\le \nu(\x-\z)= \nu(\z-\x).
\end{eqnarray*}
This shows \eqref{Lipscon}.  
\eqref{equivnorm} yields that $\dist(\cdot,C)$ is Lipschitz.  

Let $\U\subset \R^n$ be a subspace of $\R^n$ which is not contained in $C$.  Let $d(\cdot)$ be the restriction of $\dist(\cdot,C)$ to $\U$.
As $d(\cdot)$ is Lipschitz, Rademacher's theorem yields that it is differentiable almost everywhere in $\U$.  
Assume that $d(\cdot)$ is differentiable at $\x\in \U\setminus C$.
Hence $\x-\y^\star\ne \0$.  Fix $\uu\in \U$ and let $t\in\R$.  Then 
\[d(\x+t\uu)=d(\x)+\partial d(\x)(t\uu)+to(t)\le \nu(\x+t\uu-\y^\star)=\nu(\x-\y^\star)+\partial \nu(\x-\y^*)t\uu+to(t).\]
Recall that $d(\x)=\nu(\x-\y^\star)$.  Assume first that $t>0$.   Subtract $d(\x)$ from both sides of this inequality and divide by $t$.  Let $t\searrow 0$ to
deduce $\partial d(\x)(\uu)\le \partial \nu(\x-\y^*)\uu$.  By using the same arguments for $t<0$ we deduce that $\partial d(\x)(\uu)\ge \partial \nu(\x-\y^*)\uu$.
This establishes \eqref{forderdistV}.  Suppose that $\z^\star$ is another best $\nu$-approximation to $\x$.  Hence in equality \eqref{forderdistV} we can replace
$\y^\star$ by $\z^\star$.  Therefore $\partial\nu(\x-\z^\star)-\partial\nu(\x-\y^\star)$ is orthogonal to $\U$.

Assume now that $\U=\R^n$ and the norms $\nu$ and $\nu^*$ are differentiable.  Suppose that $\dist(\cdot,C)$ is differentiable at $\x\in\R^n\setminus C$.
Assume to the contrary that $\x$ has two best $\nu$-approximations $\y^\star,\z^\star$.  So $\partial\nu(\x-\y^\star)-\partial\nu(\x-\z^\star)$ is orthogonal to $\U$.
Hence $\partial\nu(\x-\y^\star)=\partial(\x-\z^\star)$.  Clearly $\dist(\x,C)=\nu(\x-\y^*)=\nu(\x-\z^\star)>0$.  Let $\x_1=\frac{1}{\nu(\x-\y^\star)}(\x-\y^\star),
\x_2=\frac{1}{\nu(\x-\z^\star)}(\x-\z^\star)\in \rS_\nu$.  So $\partial \nu(\x_1)=\partial\nu(\x_2):=\w$.  Recall that $\w\in\rS_{\nu^*}$ and hyperplanes 
$\x_i\trans \vv=1, \vv\in\R^n$ are supporting hyperplanes of $\rB_{\nu^*}$ at $\w$.  The assumption that $\nu^*$ is differentiable yields that $\x_1=\x_2$.
So $\x-\y^\star=\x-\z^\star$ which contradicts the assumption that $\y^\star\ne \z^\star$.  \qed

\section{Semi-algebraic sets}\label{sec:semalg}

Recall that a set $S\subset \R^n$ is called semi-algebraic if it is a finite union of basic semi-algebraic set given by a finite number of polynomial 
equalities $p_i(\x)=0, \ i\in \{1,...,\lambda\}$ and polynomial inequalities $q_j(\x)>0, \ j \in \{1,...,\lambda'\}$.  
The fundamental result about semi-algebraic sets (Tarski–Seidenberg theorem) is that a semi-algebraic set $S\subset\R^n$
can be described by a quantifier-free first order formula (with parameters in $\mathbb{R}$ considered as a real closed field). It follows that the
 projection of a semi-algebraic set is semi-algebraic. The class of semi-algebraic sets is closed under finite unions, finite intersections and complements.\\

A function $f:\R^n \to \R$ is called semi-algebraic if its graph $G(f)=\{(\x,f(\x)): \x\in\R^n\}$ is semi-algebraic.  We call a norm $\nu$ semi-algebraic 
if the function $\nu(\cdot)$ is semi-algebraic.  Clearly, the  norm $\|(x_1,\ldots,x_n)\trans\|_a:=
(\sum_{i=1}^n |x_i|^a)^{\frac{1}{a}}, a\ge 1$ is semi-algebraic if $a$ is rational.  Indeed, assume that $a=\frac{b}{c}$ where 
$b\ge c\ge 1$ are coprime integers, Then
\[
G(\|\cdot\|_a)=\{(x_1,\ldots,x_n,t)\trans\! : x_i=\pm y_i^c, y_i\ge 0, i=1,...,n,\  t=s^c,s\ge 0,\sum_{i=1}^n y_i^b-s^b=0\}.
\]

The following result is well known in the case when $\nu$ is the Euclidean norm \cite[\S1.1]{Cos05}. We will sketch the proof in the general case.

\begin{lemma}\label{distsemialg}  Let $C\subset \R^n$ be a nonempty closed semi-algebraic set and let $\nu$ be a semi-algebraic function. Then the function $f(\cdot):=\dist_\nu(\cdot,C)$
is semi-algebraic.
\end{lemma}
\proof  The graph of $f$ is characterized as 
\[
\{(\x,t) \in \mathbb{R}^{n+1}: t\geq 0, \  \forall \y \in C\  t \leq \nu(\x-\y), \ \forall \varepsilon > 0 \ \exists \y^\star\in C: t+\varepsilon > \nu(\x-\y^\star)\}.
\]
 This is a finite intersection of semi-algebraic sets. For example, the set $\{(\x,t) \in \mathbb{R}^{n+1}:  \forall \y \in C \  t \leq \nu(\x-\y)\}$ is 
the complement in $\mathbb{R}^{n+1}$ of the set which is the projection onto $\mathbb{R}^{n+1}$ of the set $B \subset \mathbb{R}^{n+1}\times \mathbb{R}^n$,
\[
B= (\mathbb{R}^{n+1}\times C) \cap u^{-1}(-\infty,0),
\]
where the function $u:\mathbb{R}^{n+1}\times \mathbb{R}^n\mapsto \mathbb{R}, \quad u(\x,t,\y)=\nu(\x-\y)-t$ is semi-algebraic. 
Since preimages of semi-algebraic sets by semi-algebraic maps are semi-algebraic, $B$ is semi-algebraic, and so is its projection.
A similar argument applies to  other sets in the intersection characterizing  the graph of $f$. 
\qed

For our next theorem, we will need the following results, proved as parts  of Theorem 3.3 (Fact 1 and Fact 2) in \cite{Du83}. Recall that a semi-algebraic set is called smooth
if it is an open subset of the set of smooth points of some algebraic set.

\begin{prop}\label{prop:strata} Let $X \subset \mathbb{R}^n$ be a semi-algebraic set and let $f: X \mapsto \mathbb{R}^p$ be a semi-algebraic map. Then there is a
 Whitney semi-algebraic stratification $X=\bigcup \Delta_i$ 
such that the graph $f\mid \Delta_i$ is a smooth semi-algebraic set for each $i$. 
\end{prop}

\begin{prop}\label{prop:nondiff} Let $X \subset \mathbb{R}^n$ be a smooth semi-algebraic set and let $f: X \mapsto \mathbb{R}^p$ be a map whose graph is 
 a smooth semi-algebraic set. Then the set of points in $X$ where $f$ is not differentiable is contained in a closed semi-algebraic set with dimension
 less than the dimension of $X$.
\end{prop}

Now we can prove an approximation result.  

 \begin{theo}\label{apcon}  Let $C\subset\R^n$ be a semi-algebraic set.  Assume that $\nu$ is a semi-algebraic norm such that
$\nu$ and $\nu^*$ are differentiable.
Then the set of all points $\x\in\R^n\setminus C$, denoted by $S(C)$, at which the $\nu$-approximation 
 to $\x$ in $C$ is not unique is a semi-algebraic set which does not contain an open set.
In particular $S(C)$ is contained in some hypersurface $H\subset \R^n$.
 \end{theo}

\proof
 Let $f(\x)=\dist_\nu(\x,C)$.  Since $f$ is semi-algebraic, the graph of $G(f)$ is a semi-algebraic set.  Hence 
\begin{equation}\label{defCxGf}
G(f)\times C:=\{(\x\trans,t,\y\trans)\in\R^{2n+1},\; \x\in\R^n, t=\dist_\nu(\x,C),\y\in C\}
\end{equation}
is semi-algebrac.  Let 
\begin{eqnarray}\notag
T(f):=&&\{(\x\trans,t,\y\trans,\x\trans,t,\z\trans)\trans\in\R^{2(2n+1)},
 (\x\trans,t,\y\trans,)\trans,(\x\trans,t,\z\trans,)\in G(f)\times C,\\
&&\nu(\x-\y)=\nu(\x-\z)=t,\;\|\y-\z\|^2>0\}.\label{defSing1}
\end{eqnarray}
Clearly, $T(f)$ is semi-algebraic. It is straightforward to see that $S(C)$ is the projection of $T(f)$ on the first $n$ coordinates, so  $S(C)$ is semi-algebraic.  
Theorem \ref{bestnraprxlemma}  along with Propositions \ref{prop:strata} and \ref{prop:nondiff} yields that $S(C)$  does not contain an open set.  
Hence each basic set of $S(C)$ is contained in a hypersurface.  Therefore $S(C)$ is contained
in a finite union of hypersurfaces which is a hypersurface $H$.  \qed  

\section{The case of an irreducible variety}\label{sec:irvar}
Let $p\in \F[\F^n]$.  Denote by $Z(p)$ the zero set of $p$ in $\F^n$.  
Recall that $V\subset \F^n$ is a variety (an algebraic set) if there exists a finite number of polynomials $p_1,\ldots,p_m\in \F[\F^n]$
so that $V=\cap_{i=1}^m Z(p_i)$.  Assume that $\F$ is $\R$ or $\C$ and $V\subset\F^n$ is a variety.  $V$ is called reducible if $V=V_1\cup V_2$, 
where  $V_1,V_2$ are strict subvarieties of $V$.  Otherwise $V$ is called irreducible.  Any variety $V$ has a unique decomposition as a finite union of 
irreducible varieties.  Assume that $V$ is irreducible.
Consider the Jacobian
$D(V)(\y)=([\frac{\partial p_i}{\partial x_j}(\y)]_{i=j=1}^{m,n})\trans\in \C^{n\times m}$, where $\x=(x_1,\ldots,x_n)\trans$ and $\y\in V$. 
Recall that there exist a strict subvariety of $V$, denoted $\Sing V$, such that the following conditions hold: For each $\y\in V\setminus \Sing V$
the rank of $D(V)(\y)$ is $n-d$.   Here $\Sing V:=\{\y\in V, \;\rank D(V)(\y)<n-d\}$.  The number
$d$ is the dimension of $V$ and denoted by $\dim V$.  For each $\y\in V\setminus\Sing V$ denote by $\range D(V)(\y)$ the column space of $D(\y)$.
$V$ is called smooth if $\Sing V=\emptyset$.   

Assume that $\F=\C$.
Then $V\setminus\Sing V$ is a connected complex manifold of complex dimension $d$.   
On $\C^n$ consider the symmetric form $(\x,\y):=\y\trans \x$ for $\x,\y\in\C^n$.  
The orthogonal complement of $\range D(V)(\y)$ with respect 
to $(\cdot,\cdot)$, denoted as $(\range D(V)(\y))^\perp$, is the tangent space of $V$ at $\y$.  Note that 
\[\dim \range D(V)(\y)=n-d, \quad \dim (\range D(V)(\y))^\perp= d.\] 
Let
\begin{equation}\label{defVs}
V_s:=\{\x\in V\setminus\Sing V,\; \range D(V)(\y)\cap ( \range D(V)(\y))^\perp\ne \{\0\}\}.
\end{equation} 

Let  $C\subset \R^n$ be a variety.  So $C=\cap_{i=1}^m Z(p_i)$, where $p_i\in \R[\R^n], i=1,\ldots,m$. 
Denote $C_\C:=\{\z\in \C^n,\;p_i(\z)=0, i=1,\ldots,m\}$.  In this section we assume that $C$ is irreducible.
Since $(\cdot,\cdot)$ is an inner product on $\R^n$ we deduce that
\begin{equation}\label{Vsemptyinter}
(C_\C)_s\cap (C\setminus\Sing C)=\emptyset.
\end{equation}

Fix an $\x\in\R^n$ and consider the function $g_{\x}(\y):=\|\x-\y\|^2$ restricted to $C$.  We will now study the critical points of $g_{\x}$. 
\begin{lemma}\label{critcon}  Let $C=\cap_{i=1}^m Z(p_i)\subset \R^n$  be an irreducible smooth variety.   Fix an $\x\in\R^n$.
Then $\y$ is a critical point of $g_{\x}=\|\x-\y\|^2$ on $C$ if and only if it satisfies the condition:
All minors of order $n-d+1$ of the matrix  $[D(C)(\y),\x-\y]\in\R^{n\times (m+1)}$  are zero.
\end{lemma}
\proof  Note that $\y$ is a critical point of $g_{\x}$ if an only if $\x-\y$ is orthogonal to the tangent
space of $C$ at $\y$.  I.e., $\y$ is a critical point if and only if $\x-\y\in\range D(C)(\y)$.   Since $\rank D(C)(\y)=n-d$ we deduce the lemma.
  \qed

For a general irreducible  complex variety $V\subset \C^n$ of complex dimension $d<n$ we can define the following: 
 
\begin{defn}\label{defn:semicrit} Let $V=\cap_{i=1}^m Z(p_i)\subset \C^n$ be an irreducible variety of dimension $d\in\{1,\ldots,n-1\}$.  
For each $\x\in\C^n$
we say that $\y\in V$ is a  semi-critical point corresponding to $\x$ if all minors of order $n-d+1$ of the matrix
$[D(V)(\y),\x-\y]$ are zero.
\end{defn}

Observe that a singular point of $V$ is a semi-critical point of $g_{\x}$.  
Denote by $\Sigma_1(V)\subset \C^n\times V\subset \C^{2n}$ the variety of the set of pairs $(\x,\y)\in \C^n\times V$ where  $\x$ and $\y$ satisfy 
the conditions of Definition \ref{defn:semicrit}.  For $k\in\{1,\ldots,N\}$ denote by $\pi_k: \C^N\to \C^k$
the projection of $\C^N$ onto the first $k$ coordinates.

\begin{lemma}\label{structSigmaC}  Let $V$ be an irreducible complex variety.  For $\y\in V$ denote by $\range D(V)(\y)\subset \C^n$ 
the column subspace of the matrix
$D(V)(\y)$.   Let $\Sigma_0(V)\subset \C^n\times V$ be the set all pairs $(\x,\y)\in  \C^n\times(V\setminus \Sing V)$,
where $\x\in\y+\range D(V)(\y)$.  Then $\Sigma_0(V)$ is a quasi-algebraic subset of the irreducible variety  $\Sigma(V)\subset \C^n\times V\subset\C^{2n}$
of dimension $n$, which is an irreducible component of the variety $\Sigma_1(V)\subset \C^n\times V\subset\C^{2n}$.  Each point in  $\Sigma_0(V)$ 
is a smooth point of $\Sigma(V)$. If $V$ is singular then all other irreducible components of $\Sigma(V)$ are the irreducible components of $\C^n\times
\Sing V$.
\end{lemma}

\proof  Recall that a quasi-algebraic set is a set of the form $A\setminus B$, where $A$ and $B$ are algebraic. Clearly $\C^n\times \Sing V$ 
is a subvariety $\Sigma_1(V)$.   Furthermore, $\Sigma_1(V)\setminus \C^n\times \Sing V=\Sigma_0(V)$.  
So $\Sigma_0(V)$ is  quasi-algebraic.
Since $V\setminus \Sing V$ is a manifold of dimension $d$ it follows that $\Sigma_0(V)$ is a manifold of dimension $n$.
Hence the Zariski closure of $\Sigma_0(V)$ , which is equal to the closure of $\Sigma_0(V)$  in the standard topology in $\C^n$,  is an irreducible
variety $\Sigma(V)$.  The dimension of $\Sigma(V)$ is $n$ and each point in $\Sigma_0(V)$ is a smooth point of $\Sigma(V)$.
Observe next that if $V$ is not smooth then the dimension of each irreducible component of 
$\C^n\times \Sing V$ has dimension at least $n$.  Hence no irreducible component of $\C^n\times \Sing V$ is a subvariety of $\Sigma(V)$.
Thus 
\begin{equation}\label{SigeqsmoothC} 
\Sigma_1(V)=\Sigma(V)=\Sigma_0(V) \textrm{ if } V \textrm{ is smooth.}
\end{equation}
If $V$ has singular points then the decomposition of $\Sigma_1(V)$ into its irreducible components consists of $\Sigma(V)$ and the irreducible
components of $\C^n\times \Sing V$.
\qed  

In the above notation we also have:

\begin{theo}\label{lemma:dominant}  Let $C\subset \R^n$ be an irreducible variety of dimension $d\in\{1,\ldots,n-1\}$.  Assume that $V=C_\C$.
 Let $V_s$ be given by \eqref{defVs}.  Then $V_s$ is a strict quasi-algebraic set of $V$.  For each $\y\in V\setminus (V_s\cup \Sing V)$
the map $\pi_n:\Sigma(V)\to\C^n$  is locally $1$-to-$1$ at $(y,y)$.  In particular, $\pi_n$ is dominant. 
Let
\begin{equation}\label{defUexc}
U:=\{\x\in\C^n, \; (\x,\y)\in\Sigma(V), \y\in V_s\cup \Sing(V)\}, \quad X:=\pi_n^{-1}(\C^n\setminus U)\cap \Sigma(V).
\end{equation}
Then $\pi_n^{-1}(\x)$ consists of exactly $\delta(V)$  distinct points
$(\x,\y_1(\x)),\ldots,(\x,\y_{\delta(V)})\in\Sigma(V)$  for each $\x\in \C^n\setminus U$,  i.e., $X$
is $\delta(V)$-covering of $\C^n\setminus U$. Hence  the field of rational functions over $\Sigma(V)$ is a $\delta(V)$-extension of
 the field of rational functions over $\C^n$.
\end{theo}

\proof
Clearly the map $\pi_n:\Sigma(V)\to\C^n$  is a polynomial map.  
Assume that $\y_0$ is a smooth point of $V$.  We show that $\pi_n$ is locally $1$-to-$1$ at $(\y_0,\y_0)$ if and only if 
 $\y_0\not\in V_s$.  The local coordinates of $V$ can be identified with the local coordinates $\z\in\C^d$
in the neighborhood of the origin.  So for $\y$ in a neighborhood of $\y_0 \in V$ one has $\y(\z)=\y_0+\uu(\z)$.  
Let $F(\y)$ be a fixed $n\times (n-d)$-submatrix  of $D(V)(\y)$ of rank $n-d$ in the neighborhood
of $\y_0$ in $V$.  So  the $\range D(V)(\y(\z))$ is $\{F(\y(\z))\w, \w\in\C^{n-d}\}$  for $\y(\z)=\y_0 +\uu(\z)$.  Note that $\uu(\z)=A\z+O(\|\z\|^2)$ and
$F(\y(\z))=B+O(\|\z\|)$. 
Then the projection of $\Sigma(V)$ in the neighborhood of $(\y_0,\y_0)$ onto the first $n$ coordinates
is $\y_0+A\z+B\w +O(\|\z\|^2+\|\w\|^2)$.    Hence the Jacobian of $\pi_n$ at $(\y_0,\y_0)$ is the matrix $[A\; B]$.  So $\pi_n$
is locally $1$-to-$1$ if and only if $\rank [A\; B]=n$, i.e.,  $\range D(V)(\y_0)\cap (\range D(V)(\y_0))^\perp=\{\0\}$.  
Thus $V_s$ is the set of all points in $V\setminus \Sing V$ which satisfy $\det [A\; B]=0$.  Hence $V_s$ is a quasi-algebraic.

\eqref{Vsemptyinter} yields that  $V_s$ is a strict quasi-algebraic subset of $V$. 
Furthermore,  $\pi_n$ is locally $1$-to-$1$ at $(\y_0,\y_0)$ for each $\y_0\in C\setminus \Sing C$.
Therefore $\pi_n$ is a dominant map.  The degree of this map is $\delta(V)$.  
So $\pi_n:X\to \C^n\setminus U$ is a covering map.  It is degree is $\delta(V)$.
Thus  $\pi_n^{-1}(\x)$ consists of exactly $\delta(V)$ distinct points
$(\x,\y_1(\x)),\ldots,(\x,\y_{\delta(V)})\in\Sigma(V)$  for each $\x\in \C^n\setminus U$.
Therefore the fields of rational functions over $\Sigma(V)$ is a finite extension of the field of rational functions over $\C^n$ of degree $\delta(V)$.
\qed

Let $C\subset \R^n$ be an irreducible variety of dimension $d\in\{1,\ldots,n-1\}$.
Recall the stratification of $C$ into a union of smooth quasi-algebraic sets.
With this stratification we associate the corresponding irreducible varieties in $\R^n$ as follows:
Denote $C_0:=C$.  We associate with $C$ the following directed tree $\cT(C)$.  The vertices of the tree are labeled 
by nonnegative integers from $0$ to $k\ge 0$.  Assume that $i<j$.
We say that a vertex $j$ is a descendant of the vertex $i$  (and that there is
a directed edge from $i$ to $j$) if $C_j$ is an irreducible component of $\text{Sing }C_i$.  
Denote by  $N(\{i\})\subseteq \{i+1,\ldots,k\}$  the set of all descendants of $i$.
Each $i\in\{0,\ldots,k\}$ corresponds to an irreducible subvariety $C_i$ of $C$.
The vertex $i$ is a leaf of $\cT(C)$, i.e., $i$ has no descendants,  if and only if $C_i$ is smooth.   
Otherwise, $\cup_{j\in N(\{i\})} C_j$ is the decomposition of $\Sing C_i$   into a union of its irreducible components. 
Then
\begin{equation}\label{stratV1}
C=\cup_{i\in\{0,\ldots,k\}} C_i\setminus \Sing C_i,
\end{equation}
is the smooth stratification of $C$.
\begin{defn}\label{stratV}   Let $C\subset\R^n$ be an irreducible variety.  Consider the stratification \eqref{stratV1}.
Let $V_i:=(C_i)_\C$ for $i=0,\ldots,k$.  Denote $V_0$ by $V$.
A point $\y\in V$ is called a critical point corresponding to $\x\in \C^n$ if $(\x,\y)\in\Sigma(V_i)$ for some $i\in\{0,\ldots,k\}$.
\end{defn}

We now explain briefly our definition.  Assume that $\y\in C\setminus\Sing C$.   Then $\y$ is a critical point corresponding to $\x\in\R^n$
if and only if $\y$ is a semi-critical point as defined in 
Definition \ref{defn:semicrit}.   Assume that $\y\in \Sing C$.  So $\y$ is in $C_i$ (one of the irreducible component of $\Sing C$).
Suppose furthermore that $\y_i\in C_i\setminus \Sing C_i$.  Then $\y$ can be viewed as a critical point of $g_\x|C$ if $\y$ is a critical
point of $g_\x|C_i$.
\begin{theo}\label{fibratlem}   Let $C\subset\R^n$ be a real irreducible variety. 
 Assume that the smooth stratification of $C$ is given by \eqref{stratV1}. Let $V_i:=(C_i)_\C$ for $i=0,\ldots,k$.
Let $U_i\subset \C^n$ be the variety defined in Theorem \ref{lemma:dominant} for $V=V_i$.  Let $W:=\cup_{i\in\{0,\ldots,k\}} U_i$.
For each $\x\in\R^n\setminus W$ let $Y(\x)$ be the projection of $\cup_{i\in\{0,\ldots,k\}}\pi_n^{-1}(\x)\cap( \R^{n}\times C)$
onto the last $n$ coordinates.  Then $Y(\x)$ is a nonempty finite set of $C$ of cardinality at most $\sum_{i=0}^k \delta(V_i)$.  Furthermore 
\begin{equation}\label{distforcirtpts}
\dist(\x,C)=\min_{\y\in Y(\x)} \|\x-\y\|.
\end{equation}
\end{theo}
\proof  
Recall that the projection of $\pi_n^{-1}(\z)\cap \Sigma(V_i)$ onto the last $n$ coordinates is a set of cardinality $\delta(V_i)$ in $\C^n$
for $\z\in\C^n\setminus U_i$.  Hence $|\pi_n^{-1}(\z)\cap \Sigma(V_i)|=\delta(V_i)$ for $\z\in \C^n\setminus U_i$ for $i\in\{0,\ldots,k\}$.

Let $\x\in\R^n\setminus W$.  Then $|Y(\x)|\le\sum_{i=0}^k \delta(V_i)$.  
Assume that $\dist(\x,C)=\|\x-\vv\|$ for some $\vv\in C$.  Clearly, $\vv\in C_i\setminus \Sing C_i$ for some $i\in\{0,\ldots,k\}$.
Hence $\dist(\x,C)= \dist(\x,C_i)=\|\x-\vv\|$.  Therefore $\vv$ is a critical point of $\|\x-\uu\|$, where $\uu\in C_i$.  
Thus $(\x,\vv)\in \Sigma_0(V_i)\subset \Sigma(V_i)$.  So $\vv\in Y(\x)$.  Furthermore, \eqref{distforcirtpts} holds.  \qed

\begin{theo}\label{locres}  Let $C\subset \R^n$ be an irreducible variety.  Let $V=C_\C$ and assume that 
$V_s$ is defined as in Theorem \ref{lemma:dominant}.
Let $\y\in V\setminus (V_s\cup \Sing V)$.  Then there exists an open neighborhood
$O(\y)\subset \C^n$ such that for each $\x\in O(\y)$  there exists a unique point $\y(\x)$ such that $\Sigma(C_\C)\cap (\{\x\}\times O(\y))=\{\x,\y(\x)\}$.
In particular, if $\y\in C\setminus\Sing C$ and $\x\in O(\y)\cap\R^n$ then $\dist(\x,C)=\|\x-\y(\x)\|$.
\end{theo}
\proof  Theorem \ref{lemma:dominant} claims that for each $\y\in V\setminus (V_s\cup \Sing V)$ the projection $\pi_n$ is locally $1$-to-$1$ at $(\y,\y)$.
Hence there exists an open neighborhood $O(\y)\subset \C^n$ of $\y$ such that $\pi_n: (O(\y)\times V)\cap \Sigma (V)\to O(\y)$ is a diffeomorphism. 

Assume that $\y\in C\setminus \Sing C$.   Theorem \ref{lemma:dominant} yields that $\y\in V\setminus (V_s\cup\Sing V)$.  
The proof of Theorem \ref{lemma:dominant}
implies that the Jacobian of $\pi_n$ at $(\y,\y)$ is a real matrix.  Hence $\y(\x)\in C$ for each $\x\in O(\y)\cap\R^n$ .  
Choose $O(\y)$, a small enough open neighborhood of $\y$, so that $\dist(\x,C)=\|\x-\y(\x)\|$ for each $\x\in O(\y)\cap\R^n$.
 \qed

\section{$\delta(V)$ as the top Chern number}\label{sec:smtprojvar}

Assume that $V\subset\C^n$ is a homogeneous irreducible variety.   Then $V$
induces a projective variety $V_\P\subset\P(\C^n)$.   Note that $\dim V_\P=\dim V-1=d-1$.  
It is not difficult to show that $\Sigma(V)$ is also a homogeneous variety in $\C^n\times V$.
Recall that $V_\P$ is smooth if and only if $V\setminus\{0\}$ is smooth.   Assume that $V_\P$ is smooth.
We show that the number $\delta(V)$ is equal to the top Chern number of a certain vector bundle associated with $V_\P$ under suitable conditions.   
Results  of this type are discussed in \cite{CF00,FO12,DHOST}.

View $V_\P$ as $V':=V\setminus \{0\}$ where we identify $\y\in V'$ with $t\y, t\in C\setminus\{0\}$. 
Thus $[\y]\in V_\P$ is the line spanned by $\y$ in $\C^n$.
View $\C^n$ as a trivial $n$-dimension vector bundle over $V'$.
Let $\y\in V', t\in \C\setminus\{0\}$. 
For $\y\in V'$ denote by $\cT_{V',\y}\subset \C^n$ the tangent space of $V$ at $\y$.  Then $\cT_{V,'\y}=\cT_{V',t\y}$.  
Furthermore, $\y\in\cT_{V',\y}$.   In what follows we assume that we chose $\y\in V'$ as a representative for $[\y]$.  
Let $\cB$ be the direct sum of the tangent bundle of $V_\P$ with the tautological line bundle over
$\P(\C^n)$, denoted as $\cO(-1)$. Note that the rank of $\cB$ is $d$.   Denote by $\cB'$ the dual  vector bundle of linear transformations from $\cB$ to $\C$.  
Let $\psi\in\rH^0(\cB')$  be the following global section in $\cB'$: $\psi_\y \in (\cT_{V_\P}\oplus \mathcal{O}(-1))'_{[\y]}$, 
$\psi_\y(\z)=\y\trans \z,\; \z\in \cT_{V',\y}$.
Clearly, $\psi_\y=0$ if and only if $\y$ belongs to the normal bundle of $V'$ at $\y$.
In particular $\y\trans\y=0$. \\
 Let $Q_n$ be the standard hyperquadric $Q_n=\{\z \in\C^n,\; \z\trans \z=0\}$ .  We say that $V_\P$ and $(Q_n)_\P$ 
intersect transversally if for each $\y\in V'\cap Q_n'$
the intersection of the normal bundles of $V'$ and $Q_n'$ at $\y$ is $\{0\}$.

\begin{prop}\label{transvcon}  Let $V\subset \C^n$ be a homogeneous variety.  Assume that $V_\P$ is smooth.  Let $\cB'$ and $\psi\in\rH^0(\cB')$ 
be  respectively the vector bundle
and the global section  defined above.  Then $\psi$ vanishes nowhere if and and only if $V_\P$ and $(Q_n)_\P$ intersect transversally.
\end{prop}

\proof  Clearly, the normal bundle of $(Q_n)_\P$ at $[\y]$ is given by $\span (\y)$.  Assume that $\y\in V'\cap Q_n'$.  
Then the normal bundles of $V'$ and $Q_n'$ at
$\y$ intersect transversally if and only if $\y$ is not in the normal bundle of $V'$.  That is, $\psi_\y\ne 0$.  \qed

We now give two simple examples of smooth projective varieties that intersect $Q_n$ transversally.
Let $C\subset \R^n$ be a linear subspace and denote $V=C_\C\subset \C^n$.  It is straightforward to see that 
that  $V_\P$ and $(Q_n)_\P$ intersect transversally.  (Choose an orthonormal basis in $C$.)  Consider next a quadric of the form 
$K:=\{\sum_{i=1}^n a_ix_i^2=0, \x=(x_1,\ldots,x_n)\in \C^n\}$, where $a_i\ne 0, i=1,\ldots,n$. Then $K_\P$ and $(Q_n)_\P$ intersect 
transversally if and only if  $a_i\ne a_j$ for $i\ne j$.

In the remaining part of this section we assume that $V_\P$ and $(Q_n)_\P$ intersect transversally.
For each $\x\in\C^n$ we define a global section $\phi_\x\in\rH^0(\cB')$ as follows: $\phi_{\x}(\y)$ is given by 
$\phi_{\x}(\y)(\z)=\x\trans \z,\;\z\in \cT_{V',\y}$.
Let $\bigwedge^2 \cB'$ be the the exterior product of order $2$ of the vector bundle $\cB'$ with itself.
   That is,  $\bigwedge^2 \cB'_{[\y]}$ is spanned by all $\f\wedge \g$ for  $\f,\g\in \cB'_{[\y]}$.  
Note that  $\rank \bigwedge^2 \cB'={d \choose 2}$.   Fix the section $\psi \in \rH^0(\cB')$ that was defined above. 
Denote by $\cC\subset  \bigwedge_2 \cB$ 
the subbundle induced by all global sections
$F_\x:=\phi_\x\wedge\psi\in\rH^0( \bigwedge_2 \cB')$.   Since $\psi$  vanishes nowhere it follows that $\cC$ has rank $\dim V_\P=d-1$.
We can view $F_\x$ as a global section in $\cC$.   Let $\W:=\{F_\x, \; \x\in\C^n\}$. 
For each $\y\in V'$ let $\W_\y$ be the subspace  generated by $F_{\x,\y}:=F_{\x}\mid_{[\y]}$.  Then $\W_\y=\cC_{[\y]}$.

Recall that we can associate with $\cC$ the corresponding Chern classes $c_i(\cC), j=1,\ldots,d-1$ \cite{Ha77}. 
Since $\rank \cC=\dim V_\P$ it follows that the top Chern class of $\cC$, $c_{\dim V_\P}(\cC)$, is of the from $\gamma(\cC)\omega$.
Here $\omega\in\rH^{2\dim V_\P}(V_\P,\Z)$ is the volume form on $V_\P$ such that
$\omega$ is a generator of $\rH^{2\dim V_\P}(V_\P,\Z)$ and $\gamma(\cC)\in\Z$.  
The number $\gamma(\cC)$ is called the top Chern number of $\cC$.
\begin{theo}\label{deltafor}  Let $C\subset\R^n$ be an irreducible variety.  Denote $V=C_\C$ and assume that $V$ is a homogeneus variety such that $V_\P$ is smooth.
Let $\delta(V)$ be the number of nonzero semi-critical points for $\x\in \C^n\setminus U$,  (see Definition \ref{defn:semicrit} and Theorem \ref{lemma:dominant}.)
Assume that $V_\P$ and $(Q_n)_\P$ intersect transversally.   Let $\cC$ be the vector bundle over $V_\P$ of rank $\dim V_\P$ defined above.
Then $\delta(V)$ equals to the top Chern number of $\cC$:

\begin{equation}\label{deltafor1}
\delta(V)=\gamma(\cC).
\end{equation}
\end{theo}

\proof Our theorem follows from ``Bertini type" theorem or 
Generic Smoothness Theorem \cite[Corol. III 10.7]{Ha77} and \cite[Example 3.2.16]{Ful}.  See \cite[\S2.5]{FO12} for more details.
Let $\W\subset \rH^0(\cC)$ be defined as above.  Since $\rank \cC=\dim V_\P$ and $\W$ generates $\cC$,  \cite[Theorem 2]{FO12} yields that
a generic $F_\x\in\W$ has $\gamma(\cC)$ zeros.  We now show that a generic $F_\x$ has $\delta(V)$ zeros.

Suppose first that $C$ is a subspace.  As we pointed out above, $V_\P$ and $(Q_n)_\P$ intersects transversally.
Let $T:\R^n \to C$ be the orthogonal projection on $C$.  Then $\Sigma(V)=\{(\x,T\x), \x\in\C^n\}$, i.e., $\delta(V)=1$. 
Let $\tilde V:=\{\x\in\C^n, \x\trans \y=0,\;\forall \y\in V\}$ be the orthogonal complement of $V$.   Let $\x\not\in \tilde V,\y\in V'$.
Suppose that $F_{\x,\y}=0$.  This means that the linear functionals $\phi_\x,\phi_\y:V\to C$ are linearly 
dependent.  Since $\phi_\y$ is a nonzero linear functional it follows that $\x-a\y$ is orthogonal to $V$.  Hence $T(\x)=a\y$.  
Note that $a\ne 0$ as $\x\not\in\tilde V$.
Vice versa, if $T(\x)=a\y$ then
$F_{\x}$ vanishes at $[\y]$.  Hence $F_{\x,\y}$ vanishes only at $[\y]$.  

Assume now that $C$ is not a subspace.  Then $0\in \Sing V$.  
Let $(V_\P)^\vee \subset P(\C^n)$ be the dual variety of $V_\P$ \cite{Hol88}.  That is, $[\z]\in \P(\C^n)$ is in $(V_\P)^\vee$ if and only
 if there is a point $\y\in V'$ such that $\z\in\range D(V)(\y)$.  Denote by $V^\vee\in \C^n$ the homogeneous variety induced by $(V_\P)^\vee$.
(So $(V^\vee)_\P=(V_\P)^\vee$.)
Let $U$ be the variety defined in Theorem \ref{lemma:dominant}.
Suppose that $\x\not\in U\cup V\cup V^\vee$.  Let $\y_i(\x)\in V, i=1,\ldots,\delta(V)$ be defined as in Theorem \ref{lemma:dominant}.
Since  $\C^n\times\{0\}$ is an irreducible component of $\Sigma_0(V)$ of dimension $n$, (see Lemma \ref{structSigmaC}),  
we may assume that $\y_i(\x)\ne 0$ for
$i=1,\ldots,\delta(V)$.  By definition $\x-\y_i(\x)\in \range D(V)(\y_i(\x))$.  Hence $\phi-\psi$ vanish at $\y_i(\x)$.  Therefore 
$F_\x$ vanishes at each $\y_i(\x)$.  Vice versa, suppose that $F_{\x}$ vanishes at $\y\in V'$.   
Since $\psi$ vanishes nowhere it follows that $\phi_\x-a\psi$ vanishes at $\y$.  Suppose that $a=0$.
Then $\x\in \range D(V)(\y)$, i.e., $[\x]\in (V_\P)^\vee$, contrary to our assumption. So $a\ne 0$.  Hence $a\y\in\{\y_1(x),\ldots,\y_{\delta(V)}(\x)\}$.
So $F_\x$ vanishes exactly at $[\y_i(\x)], i=1,\ldots,\delta(V)$.  \qed

\section{The distance function from an invariant subspace}\label{sec:distsubs}

\begin{theo}\label{uniqsubapprox}  Let $C\subset \R^n$ be an irreducible variety.
Let $\cG\subset \R^{n\times n}$ be a finite nontrivial  group of orthogonal matrices, $|\cG|>1$, which acts as a finite group of linear 
transformations on $\C^n$, i.e., $\x\mapsto A\x$ for $A\in\cG$.
Assume that the following conditions are satisfied:
\begin{enumerate}
\item $\U\subset \C^n$ and $\U_\R \subset \R^n$ are the subspaces of all complex respectively real fixed points of $\cG$ of complex  
respectively real dimension $\ell\in \{1,\ldots,n-1\}$; 
\item $C$ is fixed by $\cG$;
\item $\U_\R\cap C$ is a strict subset of $\U_\R$ and $C$.
\end{enumerate}
Denote $V=C_\C, C_\U:=V\cap \U$.  Then $V$ is fixed by $\cG$ and $\U\cap V$ are strict subsets of $\U$ and $V$.

Let  $V_s\subset \C^n$ be the quasi-algebraic set given by \eqref{defVs}.
Let $\y\in C_\U\setminus (V_s\cup \Sing V)$.  Then there exists a neighborhood $O(\y)\subset \C^n$ of $\y$ such that for each $\x\in O(\y)\cap \U$,
the point $\y(\x)$ given by Theorem \ref{locres} is in $C_\U$. 

Denote $\Phi:=(\U\times V)\cap \Sigma(V)$.
Then $\Phi$ contains  a finite number of irreducible components $\Phi_1,\ldots \Phi_k$, each of dimension $\ell$,
such that the the following conditions hold.    Let $\y\in C_\U\setminus (V_s\cup \Sing V)$.
Let $O(\y)\subset \C^n$ be a neighborhood of $\y$ defined as above.  Then for each $\x\in O(\y)\cap \U$ and $\y(\x)\in O(\x)\cap C_\U$ the point
$(\x,\y(\x))$ is in some $\Phi_i$.  Furthermore, for each $i\in\{1,\ldots,k\}$ there exists $\y\in C_\U\setminus (V_s\cup \Sing V)$ such that for each 
$\x\in O(\y)\cap \U$ and $\y(\x)\in O(\x)\cap C_\U$ the point $(\x,\y(\x))$ is in $\Phi_i$.  Hence each $\Phi_i$ lies in $\U\times C_\U$.
\end{theo}

\proof  Since $C$ is invariant under the action of $\cG$ it follows that $V$ is invariant under the action of $\cG$. As $\U_\R\cap C$ is a strict subset
of $\U_\R$ and $C$, it follows that $\U\cap V$ is a strict subset of $\U$ and $V$.
Assume that the group $\cG$ acts on $\C^n\times \C^n$ as follows: $(\x,\y)\mapsto (A\x,A\y)$.
Since $\cG$ is a group of orthogonal matrices, it follows that $\Sigma_0(V)\cap (\R^n\times C)$ is invariant under the action of $\cG$. 
Hence $\Sigma(V)$ is invariant under the action of $\cG$. 
So $(\x,\z)\in \Sigma(V)\iff (A\x,A\z)\in \Sigma(V)$ for each $A\in\cG$. 
Assume that $\y\in C_\U\setminus (V_s\cup \Sing V)$.   Without loss of generality we can assume that the neighborhood $O(\y)\subset \C^n$ given in  
Theorem \ref{locres} is invariant under the action of $\cG$.
 Suppose that $\x\in O(\y)\cap \U$.   Then $(\x,\y(\x))\in (O(\y)\times O(\y))\cap \Sigma(V)$.  Clearly 
$(A\x,A\y(\x))=(\x,A\y(\x))\in  (O(\y)\times O(\y))\cap \Sigma(V)$.  Theorem   \ref{locres} yields that $A\y(\x)=\y(\x)$.  Hence $\y(\x)\in C_\U$.  

For $\x\in O(\y)\cap\U$ the point $(\x,\y(\x))$ is an isolated point of $\Phi$.  Hence the set of all points $(\x,\y(\x))$, where 
$\x\in O(\y)\cap \U$, lies in $\U\times C_\U$.  
Furthermore, these points belong to an irreducible component of $\Phi$, which is denoted by $\Phi_i$.  The dimension 
of this variety is $\ell$.  Since this variety intersects $\U\times C_\U$ in an open set, it follows that $\Phi_i\subset \U\times C_\U$. 
Let $\Phi_1,\ldots,\Phi_k$ be all the irreducible components of $\Phi$ obtained in this way.   \qed

Assume that the assumptions  of Theorem \ref{uniqsubapprox} hold.     
An interesting and nontrivial problem is to give conditions such that for each $\x\in\U_\R$  there exists a best $C$-approximation 
in $C\cap \U_\R$.  
We now give sufficient conditions so  that this property holds.  Recall that a set $Q$ in $\R^n$ or $\C^n$ is called constructible if it is a finite union 
of quasi-algebraic sets.  Chevalley's theorem states that the image of a constructible set by a polynomial map is a constructible set.

\begin{theo}\label{sufconglob}  Let the assumptions of Theorem \ref{uniqsubapprox} hold.    Let $\Phi_1,\ldots,\Phi_k$ be the irreducible
components of $\Phi$ defined in Theorem \ref{uniqsubapprox}.  Let $\Phi=\Phi'\cup(\cup_{i=1}^k\Phi_i)$ 
be the decomposition of $\Phi$ into its irreducible components.  Assume the following conditions.  First, the equality
\begin{equation}
\dist(\x,C)=\dist(\x,\Sing C) \textrm{ for } \x\in\U_\R
\end{equation}
defines a semi-algebraic $Q$ of dimension less than $\ell$ in $\U_\R$.  ($Q=\emptyset$ if $\Sing C=\emptyset$.)  
Second, the constructible set $\pi_n(\Phi'\cap(\R^n\times\R^n))$ has 
dimension less than $\ell$, i.e., the closure of $\pi_n(\Phi'\cap(\R^n\times \R^n))$ in $\U_\R$ is a strict subvariety of $\U_\R$.  
Then for each $\x\in\U_\R$ there exists a $C$-best approximation which is in $C\cap \U_\R$.
\end{theo} 
\proof   Clearly  $Q$ is a closed semi-algebraic set.
Since $\U_\R\setminus Q$ is an open set which is dense in $\U_\R$ it follows that for each $\x\in\U_\R$
\begin{equation}\label{formdistU}
\dist(\x,C)=\min\{\|\x-\y\|,\; \x\in\U_\R,  (\x,\y)\in\Phi\cap (\R^n\times \R^n)\}.
\end{equation}
The assumption on $\Phi'$ means that 
\[\pi_n((\Phi\setminus\Phi')\cap(\R^n\times \R^n))=\pi_n(\cup_{i=1}^k(\Phi_i\setminus\Phi')\cap(\R^n\times \R^n))\subset \U_R\times (C\cap \U_\R)\]
is a constructible set in $\U_\R$ which is dense in $\U_\R$.  Observe next that for each $\x\in \pi_n((\Phi\setminus\Phi')\cap(\R^n\times \R^n))\setminus Q$
each $C$-best approximation is in $C\cap U_\R$.  As  $\pi_n((\Phi\setminus\Phi')\cap(\R^n\times \R^n))\setminus Q$ is dense in $\U_\R$ it follows that
for each $\x\in\U_\R$ there exists a $C$-best approximation which is in $C\cap \U_\R$.  \qed

Two cases where the conclusion of Theorem \ref{sufconglob} hold are discussed in \S\ref{sec:aprrten}.
We close this section with the following simpler result:
\begin{prop}\label{sufconloc}  Let the assumptions of Theorem \ref{uniqsubapprox} hold.   Then
there exists a semi-algebraic open set $O\subset \R^n$ containing all smooth points  of $C$ lying in $C\cap\U_\R$ such that for each
$\x\in O$ a best $C$ approximation of  $\x$ is unique.  Furthermore, for each $\x\in O\cap\U_\R$ a best $C$ approximation lies in $C\cap\U_\R$ .

\end{prop}
\proof  Let $S(C)\subset \R^n$ be the set of all points for which a best $C$ approximation of  $\x$ is not unique. 
Theorem  \ref{apcon} claims that $S(C)$ is semi-algebraic, and it is contained in some hypersurface $H\subset R^n$.
Recall that a complement of a semi-algebraic set is semi-algebraic.  Let $O_1:=\R^n\setminus S(C)$.  So each point $\x\in O_1$
has a unique best $C$-approximation.  Let $O_2$ by the set of all $\x\in O_1\cap\U_\R$ such that a best $C$ approximation is in $C_{\U_\R}$.
The arguments of the proof of Theorem \ref{apcon} yield that $O_2$ is semi-algebraic.

Since $\U_\R \cap C$ is a strict subset $\U_\R$ and $C$ it follows that $\U_\R\cap (C\setminus\Sing C)\ne\emptyset$. 
Let $V_s$ be defined as in Theorem \ref{lemma:dominant}.   Theorem \ref{lemma:dominant} claims that $V_s\cap (C\setminus \Sing C)=\emptyset$.
Theorems \ref{locres} and \ref{uniqsubapprox} yield that $O_2$ contains a neighborhood of each $\y\in \U_\R\cap (C\setminus  \Sing C)$.
Let $O_3:=\R^n\setminus O_2$.  So $O_3$ is semi-algebraic.  Recall that a closure of a semi-algebraic set is semi-algebraic.
So Closure$(O_3)$ is semi-algebraic.  Hence $O:=\R^n\setminus$Closure$(O_3)$, the interior of $O_2$, is a set satisfying the conditions of the proposition.
\qed

\section{On the rank of generic  tensors in  $\F^{m^{\times d}}$}\label{sec:genrankt}

Let $\F^{m^{\times d}}\supset \rS(m,d,\F)$ be the space of $d$-mode, ($d \ge 2$), tensors whose each mode 
has $m$ coordinates, and the subspace of $d$-mode symmetric tensors over a field $\cF$.   
For $\x_1,\ldots,\x_d\in\F^m$ denote by $\otimes_{i=1}^d\x_i$ the tensor product $\x_1\otimes\ldots\otimes\x_d\in \F^{m^{\times d}}$.
If $\x_i=\uu$ for $i=1,\ldots,d$ then $\otimes^d \uu:=\otimes_{i=1}^d\x_i$.
Let $\cT\in\F^{m^{\times d}}\setminus\{0\}, \cS\in\rS(m,d,\F)\setminus\{0\}$.
Consider the following decompositions of $\cT$ and $\cS $ into rank one tensors:
\begin{eqnarray}\label{rankten}
&&\cT=\sum_{j=1}^r \otimes_{i=1}^d \x_{i,j},\quad \x_{i,j}\in\F^m\setminus\{\0\},\;i=1,\ldots,d, j=1,\ldots,r, \\
&&\label{rankrsymten}
\cS=\sum_{j=1}^s t_j \otimes^d\uu_j, \quad t_j\in\F\setminus\{0\},\uu_j\in\F^m\setminus\{\0\}, j=1,\ldots,s.
\end{eqnarray}
The minimal $r$ and $s$ for which the above equalities holds for $\cT$ are called the \emph{rank} and  \emph{symmetric rank} of $\cT$
and $\cS$ respectively,  which are denoted by $\rank \cT$ and $\srank \cS$.   (The rank and the symmetric rank of zero tensor is zero.)
Note that if $\F$ is algebraically closed that in \eqref{rankrsymten} we can assume that each $t_j=1$.  For $\F=\R$ we can assume that each $t_j=\pm 1$.

 Let $\e_i:=(\delta_{i1},\ldots,\delta_{im})\trans, i=1,\ldots,$
be the standard basis in $\F^m$.  Then $\otimes_{j=1}^d \e_{i_j}$, $i_1,\ldots,i_d=1,\ldots,m$ is a standard basis in  $\F^{m^{\times d}}$.
Hence any $\cT\in \F^{m^{\times d}}$ has a decomposition \eqref{rankten}.  So $\rank\cT$ is well defined.

However the symmetric rank of a symmetric tensor $\cS\in\rS(m,d,\F)$ may be not defined for a field of a finite characteristic $p\ge 2$.
We now show that such $\cS$ exist for any finite field $\F_{p^l}$ with $p^l$ elements and a corresponding $d$.

Denote by $\rS_d$ the group of permutations on $\{1,\ldots,d\}$.  Let  $1\le i_1\le \ldots\le i_d\le m$.  Denote by orb$(i_1,\ldots,i_d)$ the orbit
of the multiset $\{i_1,\ldots,i_d\}\subset \{1,\ldots,m\}^d$ under the action of $\rS_d$.  I.e., this orbit is a union of all ordered distinct 
multisets $\{i_{\sigma(1)},\ldots,i_{\sigma(d)}\},\sigma\in\rS_d$.  Note that the number of such orbits is ${m+d-1\choose d}$.
To each such orbit we associate a symmetric tensor $\cS(i_1,\ldots,i_d):=\sum_{\{j_1,\ldots,j_d\}\in{\rm orb}(i_1,\ldots,i_d)} \otimes_{k=1}^d \e_{j_k}$.
Clearly the set of all these symmetric tensors form a basis in $\rS(m,d,\F)$.  Hence
$\dim\rS(m,d,\F)={m+d-1\choose d}$. 
(The following two results are probably known, and we give their proof for completeness.)
\begin{prop}\label{noexsymdecomp}  Let $p\ge 2$ be a prime, $\F_{p^l}$ be a field with $p^l$ elements and
 $m\ge 2$ be an integer.  Then for each integer $d$ satisfying 
\begin{equation}\label{dineq}
{m+d-1\choose d}>\frac{p^{ml}-1}{p^l-1}
\end{equation}
there exists $\cS\in\rS(m,d,\F_{p^l})$ such that $\cS$ is not a linear combination of rank one symmetric tensors.
\end{prop}

\proof  The number of nonzero elements in $\F_{p^l}^m$ is $p^{ml}-1$.  
Each nonzero vector $\x\in \F_{p^l}^m$ generates a line with $p^{l}-1$ nonzero elements of the form $t\x,t\in\F_{p^l}\setminus\{0\}$.
Hence the subspace $\U\subset\rS(m,d,\F_{p^l})$ generated by all vectors of the form 
$\otimes^d \uu, \uu\in\F_{p^l}^m$ is of dimension $\frac{p^{ml}-1}{p^l-1}$ at most.  Assume that \eqref{dineq} holds.  So $\dim \U <\dim \rS(m,d,\F_{p^l})$.
Hence there exists $\cS\in \rS(m,d,\F_{p^l})\setminus \U$.  \qed

\begin{prop}\label{existsymrankcharF0}  Let $\F$ be a field with at least $d$ elements.  Then for each $\cS\in\rS(m,d,\F)$ \eqref{rankrsymten} holds.
\end{prop}
\proof  For $\uu,\vv\in \F^m$ and a variable $t\in\F$ let
\[\otimes^d(t\uu+\vv)=\sum_{k=0}^d t^k\cS_{k,d-k}(\uu,\vv).\]
So $\cS_{k,d-k}(\uu,\vv)$ is the symmetric tensor induced by $(\otimes^k\uu)\otimes(\otimes^{d-k}\vv)$.  Clearly
$\cS_{0,d}(\uu,\vv)=\otimes^d\vv, \cS_{d,0}(\uu,\vv)=\otimes^d\uu$.
We claim that $\cS_{k,d-k}(\uu,\vv)$ has a decomposition \eqref{rankrsymten} for $k=1,\ldots,d-1$ .  Consider the polynomial 
\[\otimes^d(t\uu+\vv) -t^d\otimes^d\uu=\sum_{k=0}^{d-1} t^k\cS_{k,d-k}(\uu,\vv).\]
Recall that the Vandermonde matrix $[\tau_i^j]_{i,j=0}^{d-1}$ is an invertible matrix for a subset  $\{\tau_0,\ldots,\tau_{d-1}\}$ of $\F$
of cardinality $d$.
Hence each $\cS_{k,d-k}(\uu,\vv), k=0,\ldots,d-1$ can be expressed as a linear combination of $\otimes^d(\tau_i\uu+\vv), i=0,\ldots,d-1$ and $\otimes^d\uu$.

More generally, consider the following polynomial 
\[\cX(t_1,\ldots,t_{d-1}):=\otimes^d(\e_m+\sum_{j=1}^{m-1}t_j\e_i)-\sum_{j=1}^{m-1}t_j^d\otimes^d\e_j\]
 in $d-1$ variables $t_1,\ldots,t_{m-1}$.   The coefficients the monomial $t_1^{j_1}\ldots t_{d-1}^{j_{d-1}}$ is the tensor 
$\cS(i_1,\ldots,i_d)$  in the standard basis of $\rS(m,d,\F)$ described above.  Vice versa, each vector in the standard basis of $\rS(m,d,\F)$, except 
$\otimes^d \e_j, j=1,\ldots,m-1$, is a coefficient of the corresponding monomial
$t_1^{j_1}\ldots t_{d-1}^{j_{d-1}}$.   View $\cX(t_1,\ldots,t_{m-1})$ as a polynomial in $t_{m-1}$:  
\[\cX(t_1,\ldots,t_{m-1})=\sum_{k=0}^{d-1} t_{m-1}^k \cX_k(t_1,\ldots,t_{m-2}).\] 
Find the polynomials $ \cX_k(t_1,\ldots,t_{m-2})$ for $k=0,\ldots,d-1$ using the above procedure.  Continue this procedure to
obtain all vectors in the standard basis of $\rS(m,d,\F)$, except $\otimes^d \e_j, j=1,\ldots,m-1$ as linear combinations of vectors
$\cX(t_1,\ldots,t_{m-1})$, where $t_1,\ldots,t_{m-1}\in\{\tau_0,\ldots,\tau_{d-1}\}$.   Then  \eqref{rankrsymten} holds. \qed

 For the case where $\F$ is an infinite field see \cite{AH95}.
\begin{corol}\label{spanrnk1symten}  Let $m,d\ge 2$ be integers.   Let $\F$ be a field with at least $d$ elements.  Then there exists $k={m+d-1\choose d}$
vectors $\uu_1,\ldots,\uu_k\in\F^m$ such that $\otimes^d\uu_1,\ldots,\otimes^d\uu_k$ form a basis in $\rS(m,d,\F)$.
\end{corol}

Assume that $\F$ is an infinite field.  Clearly 
\begin{equation}\label{rankin}
\rank \cS\le \srank \cS \textrm{ for } \cS\in \rS(m,d,\F).
\end{equation} 
\begin{example}\label{countexsymrankcon}  Let $\F$ be any field of characteristic $2$.  Let $\cS=\e_1\otimes\e_2+\e_2\otimes\e_1\in\rS(2,2,\F)$.
Then
\begin{equation}\label{countexsymrankcon1}
\rank\cS=2, \quad \srank\cS=3.
\end{equation}
\end{example}
\proof  Clearly $\rank\cS=2$.  The equality $\cS=\e_1\otimes \e_1+\e_2\otimes\e_2+(\e_1+\e_2)\otimes(\e_1+\e_2)$ yields that $\srank\cS\le 3$. 
Assume now that $\srank\cS=2$.  So $\cS=a\uu\otimes\uu+b\vv\otimes\vv$.  Since $\rank\cS=2$ we deduce that $\uu,\vv\in\F^2$ are linearly
independent and $a,b\ne 0$.  Let $\uu=(u_1,u_2)\trans,\vv=(v_1,v_2)\trans$.  Clearly 
\[a u_1^2+bv_1^2=a u_2^2+b u_2^2=0.\]
If $u_1=0$ then $v_1=0$ which contradicts the assumption that $\uu,\vv$ are linearly independent.  Hence $u_1,v_1,u_2,v_2\ne 0$.
Since the characteristic of $\F$ is $2$ it follows 
\[\frac{a}{b}=\frac{v_1^2}{u_1^2}=\frac{v_2^2}{u_2^2}\Rightarrow 0=v_1^2u_2^2+v_2^2u_2^2=(v_1u_2+v_2u_1)^2\]
So $\uu,\vv$ are linearly dependent contrary to our assumption.  \qed

It is an open problem if equality holds in \eqref{rankin} for each symmetric tensor $\cS$ for $\F=\C$, or more generally over any algebraically closed (or
an infinite) field $\F$ of characteristic different from $2$.  (This is true for $d=2$.)   \cite[Proposition 5.5]{CGLM08} 
shows that  equality holds in \eqref{rankin} for $\F=\C$ and $\srank\cT\le 2$.

\begin{defn}\label{generdef}   Let  $\F$ be a field, $m\ge 2, d\ge 3$ be integers.  Assume that  $\cT\in\F^{m^{\times d}}\setminus\{0\},
\cS\in\rS(m,d,\F)\setminus\{0\}$ have decompositions \eqref{rankten} and \eqref{rankrsymten} respectively, which can be rewritten as follows:
\begin{eqnarray}\label{rankten1}
&&\cT=\sum_{j=1}^r (\otimes_{i=1}^{a} \x_{i,j})\otimes  (\otimes_{i=a+1}^{a+b} \x_{i,j})\otimes (\otimes_{i=a+b+1}^{a+b+c} \x_{i,j}) , \\
&&\label{rankrsymten1}
\cS=\sum_{j=1}^s t_j (\otimes^{a}\uu_j)\otimes (\otimes^{b}\uu_j)\otimes(\otimes^{c}\uu_j).
\end{eqnarray}
Here $a,b,c$ are positive integers such that $a+b+c=d$.
Then $\cT$ and $\cS$ of are called $(a,b,c)$-generic if the following conditions hold for $\cT$ and $\cS$ respectively.
\begin{enumerate}
\item  For $\cT$: any $\min(m^a,r)$ vectors out of $\otimes_{i=1}^a\x_{i,1},\ldots,\otimes_{i=1}^{a}\x_{i,r}$ are linearly independent; any
$\min(m^{b},r)$ vectors out of $\otimes_{i=a+1}^{a+b}\x_{i,1},\ldots,\otimes_{i=a+1}^{a+b}\x_{i,r}$  are linearly independent; 
any $\min(m^{c },r)$ vectors out of
$\otimes_{i=a+b+1}^{d}\x_{i,1}$, $\ldots$ ,$\otimes_{i=a+b+1}^{d}\x_{i,r}$ are linearly independent.
\item For $\cS$: if $\F$ is algebraically closed then each $t_j=1$ and for $\F=\R$ each $t_j=\pm 1$;
any $\min({m+a-1\choose a}, s)$ vectors out of $\{\otimes^a\uu_1,\ldots,\otimes^a\uu_{s}\}$ are linearly independent; any $\min({m+b-1\choose b},s)$ out of 
$\otimes^b\uu_1,\ldots,\otimes^b\uu_s$ are linearly independent; any $\min({m+c-1\choose c}, s)$ vectors out of $\{\otimes^c\uu_1,\ldots,\otimes^c\uu_{s}\}$
 are linearly independent.
\end{enumerate}
\end{defn} 
 
Assume that $\F=\C$,  $d=3$ and $r\le m$.  Then the Definition \ref{generdef} of a generic symmetric tensor coincides with the definition in \cite{CGLM08}.
\begin{theo}\label{rankgenten}  Let $\F$ be a field, $m\ge 2, d\ge 3$ be integers.  Assume that the decompositions (\ref{rankten1}-\ref{rankrsymten1}) of
$\cT\in\F^{m^{\times d}}$ and $\cS\in\rS(m,d,\F)$ are $(1,\lfloor\frac{d-1}{2}\rfloor, \lceil\frac{d-1}{2}\rceil)$-generic.  
\begin{enumerate}
\item Let $d=2b+1\ge 3$.  Then
\begin{eqnarray}\label{rankgenT=r1}
\rank \cT=r \textrm{ if } r\le m^b +\frac{m-2}{2},\\
\label{srankS=s1}
\rank \cS=\srank \cS=s \textrm{ if } r\le {m+b-1\choose b}+\frac{m-2}{2}.
\end{eqnarray}
\item Let $d=2b+2\ge 4$.  Then
\begin{eqnarray}\label{rankgenT=r2}
\rank \cT=r \textrm{ if } r\le m^b +m-2,\\
\label{srankS=s1}
\rank \cS=\srank \cS=s \textrm{ if } r\le {m+b-1\choose b}+m-2.
\end{eqnarray}
\end{enumerate}
In all the above cases representation of $\cT$ and $\cS$ as a sum of rank one tensors is unique up to a permutation of the summands.
\end{theo}
\proof  The proof of the theorem uses Kruskal's theorem for $3$-mode tensors $\F^{n\times p\times q}:=\F^n\otimes\F^p\otimes\F^q$
\cite{Kru77}.  (See \cite{Rho10} for a short proof of Kruskal's theorem.)   Let $\x_1,\ldots,\x_l\in\F^n$ be nonzero vectors.  Form the matrix $X=[\x_1\ldots
\x_l]\in\F^{n\times l}$, (whose columns are $\x_1,\ldots,\x_l$).   Then Kruskal rank of $X$, denoted as $\Krank X$, is the maximal integer $k\le l$ such that
any $k$ vectors out of $\{\x_1,\ldots,\x_l\}$ are linearly independent.  Let $\cX\in\F^{n\times p\times q}$ and assume that 
\begin{equation}\label{dec3tenX}
\cX=\sum_{i=1}^r \y_i\otimes\z_i\otimes \w_i,\quad \y_i\in\F^n, \z_i\in\F^p,\w_i\in\F^q, i=1,\ldots,r.
\end{equation}
Form the matrices $ Y=[\y_1\ldots\y_r],Z=[\z_1\ldots\z_r], W=[\w_1\ldots\w_r]$.  Let $\kappa_1=\Krank Y, \kappa_2=\Krank Z, \kappa_3=\Krank W$.
We call $(\kappa_1,\kappa_2,\kappa_3)$ the Kruskal ranks of the decomposition \eqref{dec3tenX}.
Kruskal's theorem claims that if 
\begin{equation}\label{kruskalcon}
\kappa_1+\kappa_2+\kappa_3\ge 2r+2
\end{equation}
then $\rank \cX=r$.  
Furthermore, the decomposition of $\cX$ to a sum of $r$ rank one tensors is unique up to a permutation of  the summands. 

We now prove our theorem.   
Assume that the decompositions of tensors $\cT$ and $\cS$ given by  (\ref{rankten1}-\ref{rankrsymten1}) are $(a,b,c)$ generic.
(For $\cS$ we use the identity $t_j(\otimes^a\uu_j)\otimes (\otimes^b\uu_j)\otimes(\otimes^c\uu_j)=
(t_j\otimes^a\uu_j)\otimes (\otimes^b\uu_j)\otimes(\otimes^c\uu_j)$ for $j=1,\ldots,r$.)
Then the Kruskal ranks of  these decompositions of $\cT$ and $\cS$, denoted by $(\alpha_1,\alpha_2,\alpha_3)$ and $(\beta_1,\beta_2,\beta_3)$ respectively,
are:
\begin{eqnarray*}
&&\alpha_1= \min(m^a,r), \alpha_2=\min(m^b,r), \alpha_3=\min(m^c,r),\\
&&\beta_1=\min({m+a-1\choose a},r), \beta_2=\min({m+b-1\choose b},r),\beta_3=\min({m+c-1\choose c},r).
\end{eqnarray*}

For $r=1$ the theorem is trivial.  In what follows we assume that $r\ge 2$.
Suppose first that $d=2b+1$.  Assume first that $\cS$ is $(1,b,b)$ generic.   Suppose first that $r\le m$.
Then $\alpha_1=\alpha_2=\alpha_3=r\ge 2$.  So Kruskal's condition \eqref{kruskalcon} holds.  Let $\cS$ be represented as a sum of $r'$
rank one tensors in $\F^{m^{\times d}}$, where  $r'\le r$.  By viewing $\cS$ as a tensor in $\F^{m^{\times a}}\otimes \F^{m^{\times b}}\otimes \F^{m^{\times c}}$
and using Kruskal's theorem we deduce that $r'=r$.   Moreover, the decomposition of $\cS$ to a sum of rank one tensors in $\F^{m^{\times d}}$ is unique up to a permutation
of summands.

Assume now that  $m\le r\le {m+b-1\choose b}$.  Then $\alpha_1=m, \alpha_2=\alpha_3=r$.  Again, Kruskal's inequality holds. 
We deduce the theorem in this case for $\cS$.  Suppose that $m\ge 4$ and $r> {m+b-1\choose b}$.  Then $\alpha_1=m, \alpha_2=\alpha_3={m+b-1\choose b}$.
(Recall that $\dim \rS(m,d,\F)={m+b-1\choose b}$.)
Then Kruskal inequality holds if and only if $r\le {m+b-1\choose b}+\frac{m-2}{2}$.  Hence the theorem holds in this case for $\cS$ too .

Assume now that $d=2b+2$.   The above arguments apply for $\cS$ and $r\le {m+b-1\choose b}$.  Assume that $m>2$ and ${m+b-1\choose b}<r\le {m+b\choose b+1}$.
Then $\alpha_1=m, \alpha_2={m+b-1\choose b}, \alpha_3=r$.   Kruskal's inequality yields that $r\le {m+b-1\choose b}+m-2$.  Hence the theorem holds in this case too
for $\cS$.  Similar arguments yield the theorem for $\cT$.  \qed

The upper bound on $r$ in Theorem \ref{rankgenten} for $\cS$ can be replaced by  
\begin{equation}\label{defNnd}
N(m,d)= 
\begin{cases}
{m+\frac{d-3}{2}\choose m-1}+\frac{m-2}{2} \textrm{ for an odd } d\ge 3\\
{m+\frac{d-4}{2}\choose m-1}+m-2\textrm{ for an even } d\ge 4
\end{cases}
\end{equation}
(We  used the identity ${n \choose k}={n\choose n-k}$.)

 \section{Approximation of real symmetric tensors}\label{sec:aprrten}

Let $\cR(k,m^{\times d})\subset \R^{m^{\times d}}, \cR_\C(k,m^{\times d})\subset \C^{m^{\times d}}$ 
be the closure of all tensors of rank at most $k$,  i.e., all tensors of border rank at most  $k$ in $\R^{m^{\times d}},\C^{m^{\times d}}$ respectively.
In this section we consider the best approximation problem in $\R^{m^{\times d}}$ for $C:=\cR(k,m^{\times d})$.
Clearly, $C_\C= \cR_\C(k,m^{\times d})$.
Observe next that the symmetric group $\rS_d$ of order $d$ acts on $\C^{m^{\times d}}$ as follows:  For $\sigma\in \rS_d$
and $\cT=[t_{i_1,\ldots,i_d}]\in \C^{m^{\times d}}$ we define
$\sigma(\cT)=[t_{i_{\sigma(1)},\ldots,i_{\sigma(i_d)}}]$.  
Clearly, the action of each $\sigma$ preserves the Hilbert-Schmidt norm on $\R^{m^{\times d}}$.
It is straightforward to see that   $C$ and $C_\C$ are invariant under the action of $\rS_d$.
Furthermore $\rS(m,d,\F)$ is the set of the fixed points in $\F^{m^{\times d}}$ for $\F=\R,\C$ respectively.  

A natural question is if for each $\cS\in\rS(m,d,\R)$ there exists a best $k$-border rank approximation which is symmetric.
This is a special case of the approximation problem discussed in \S\ref{sec:distsubs}.  

For $d=2$, i.e., the space of real symmetric matrices,  the answer to this question
 is positive \cite{GV}. 
Let $C_k$ is the set of real matrices of rank at most $k$.  Then the set of singular points of $C_k$ is $\Sing C_k=C_{k-1}$.    
A best $k$-rank approximation of $B\in\rS(m,2,\R)$
is $A\in C_k$ which has the same $k$ maximal singular values and corresponding left and right singular vectors as $B$.  
Hence the set $Q$ defined in
Theorem \ref{sufconglob} is $C_{k-1}$.  Assume now that the $m$ eigenvalues of $B$, $\lambda_1,\ldots,\lambda_m$ 
satisfy the condition $|\lambda_i|\ne
|\lambda_j|$ for $i\ne j$.  (I.e., all singular values of $B$ are distinct.)  So 
$B=\sum_{i=1}^m \lambda_i \uu_i\uu_i\trans$, where $\uu_i\trans\uu_j=\delta_{ij},
i,j=1,\ldots,m.$
Then $A$ is a critical point of the function $\tr (B-X)\trans (B-X), X\in C_k$ if and only
if $A=\sum_{i\in\Omega} \lambda_i\uu_i\uu_i\trans$ for any subset $\Omega$ of $\{1,\ldots,m\}$ of cardinality $k$.
So each real critical point of $B$ is symmetric.  Hence the assumptions of Theorem \ref{sufconglob} hold in this case.

For $k=1$ and any $d\ge 3$ the answer to this question is also positive.  
I.e., every symmetric tensor  $\cS\in\rS(d,m,\R)$ has a symmetric best rank one approximation.
This result is implied by Banach's theorem \cite{Ban38}.
The Banach theorem was re-proved in \cite{CHLZ, Fri13}.
The results in \cite{FO12} yield that the set of symmetric tensors which do not have a unique best rank one approximation has zero Lebegue measure.
We now give an improved version of this result.
\begin{theo}\label{auniqbsr1apsymt}  Let $m\ge 2, d\ge 3$ be integers.  Then there exists a semi-algebraic set $W\subset \rS(m,d,\R)$  which does not contain
an open set  with the following properties.
$\cS\in\rS(m,d,\R)$ does not have a unique best rank one approximation if and only if $\cS\in W$.  Furthermore, for each $\cS\in \rS(m,d,\R)\setminus W$
the unique best rank one approximation is symmetric.  In particular, there exists a hypersurface $H\subset \rS(m,d,\R)$ such that 
 for each $\cS\in \rS(m,d,\R)\setminus H$
the unique best rank one approximation is symmetric.
\end{theo}
\proof  
Let $W\subset \rS(m,d,\R)$ be the set of all $\cS\in\rS(m,d,\R)$ that do not have a unique best rank one approximation.
The arguments of the proof of Theorem \ref{apcon} imply that $W$ is semi-algebraic.  
Let $\cS\in \rS(m,d,\R)\setminus W$.  Then $\cS$ has a unique best rank one approximation.
\cite{Fri13} claims that $\cS$ has a best rank one approximation which is symmetric.   Hence $\cS$ has a unique best rank one approximation
which is symmetric.  

We now show that $W$ does not contain an open set.
Recall that the function $\dist(\cdot, \cR(1,m^{\times d})):\R^{m^{\times d}}\to \R$ is semi-algebraic.
Let $d(\cdot):\rS(m,d,\R)\to \R$ be the restriction of  $\dist(\cdot, \cR(1,m^{\times d}))$ to $\rS(m,d,\R)$.
Clearly, $d(\cdot)$ is semi-algebraic.  Proposition \ref{prop:nondiff} yields that $d(\cdot)$ is not differentiable on a semi-algebraic set
$W'\subset \rS(m,d,\R)$ of dimension less than $\dim \rS(m,d,\R)$.  

Let $\cS\in \rS(m,d,\R)\setminus W'$.  Hence function $d(\cdot)$ is differentiable at $\cS$.  The arguments of \cite[\S7]{FO12} yield that
 the set of all best rank one approximation of $\cS$
is the orbit of one best rank one approximation $\otimes_{i\in[d]}\x_i$  under the action of the symmetric group $\rS_d$ as defined above.  
That is, all best rank one approximations are of the form  $\otimes_{i\in[d]}\x_{\sigma(i)}$ for $\sigma\in\rS_d$.
\cite{Fri13} claims that 
$\cS$ has a best rank one approximation which is symmetric.   Hence the orbit of $\otimes_{i\in[d]}\x_i$ consists of one symmetric tensor.
In particular, $\cS$ has a unique  best rank one approximation which is symmetric.  Therefore $W'\supset W$.  Hence $W$ does not contain an open set.
The arguments of the proof of Theorem \ref{apcon} imply that $W$ is contained in a hypersurface $H$.
\qed

For rank one approximation of tensors in $\rS(m,d,\R)$ the assumptions of Theorem \ref{sufconglob} are equivalent to the following statement.
There exists an open quasi-algebraic set $W\subset \rS(m,d,\R)$ such that for each $\cB\in W$  every real rank one tensor $\cA$ which is a critical point of $\|\cB-\cX\|^2,
\cX\in\cR(1,m^{\times d})$ is symmetric. 

For $k\ge 2$ and $d\ge 3$ we have the following weaker result.
\begin{theo}\label{brraproxsymt}  Let $m\ge 2, d\ge 3$ be integers.  Assume that $2\le k\le N(m,d)$.  Then there exists an open semi-algebraic set $O\subset \rS(m,d)$
containing all $(1,\lfloor\frac{d-1}{2}\rfloor,\lceil\frac{d-1}{2}\rceil)$-generic real symmetric tensors of rank $k$ such that for
each $\cS\in O$ a best $k$-border rank approximation is symmetric and unique.  

\end{theo}
\proof  Let $\cS$ be an $(a,b,c)$ symmetric generic tensor given by \eqref{rankrsymten1}, where $a=1,b=\lfloor\frac{d-1}{2}\rfloor, 
c=\lceil\frac{d-1}{2}\rceil$, $r=k\in [2,N(m,d)]$ and each $t_j=\pm 1$.    We claim that $\cS$ is a smooth point of $\cR(k,m^{\times d})$.
Recall the proof of Theorem \eqref{rankgenten}. 
Kruskal's theorem yields that any $\cT\in\cR(k,m^{\times d})$ in a suitably small neighborhood of $\cS$ must be of the form \eqref{rankten1}, such that each
$(\otimes_{i=1}^{a} \x_{i,j})\otimes  (\otimes_{i=a+1}^{a+b} \x_{i,j})\otimes (\otimes_{i=a+b+1}^{a+b+c} \x_{i,j})$  is in the neighborhood
of $t_j (\otimes^{a}\uu_j)\otimes (\otimes^{b}\uu_j)\otimes(\otimes^{c}\uu_j)$.  Hence $\cT$ is $(a,b,c)$-generic and its rank is $k$.
Therefore $\cS$ is a smooth point of $\cR(k,m^{\times d})$.  Theorem  \ref{locres} and Proposition \ref{sufconloc} imply our theorem.  \qed

\end{document}